\documentclass[reqno,11pt]{amsart} 
\usepackage{amsmath}
\usepackage{amssymb}
\usepackage{amsthm}
\usepackage{verbatim}
\usepackage{xcolor}
\usepackage{graphics}

\newcommand\NoBlackBoxes{\global\overfullrule0pt}

\NoBlackBoxes
\theoremstyle{plain}

\begin{document}

\title{QUANTIFIED CRAM\'ER-WOLD CONTINUITY \\
THEOREM FOR THE KANTOROVICH \\ 
TRANSPORT DISTANCE
}

\author{Sergey G. Bobkov$^{1,3}$}
\thanks{1) 
School of Mathematics, University of Minnesota, Minneapolis, MN, USA,
bobkov@math.umn.edu. 
}

\author{Friedrich G\"otze$^{2,3}$}
\thanks{2) Faculty of Mathematics,
Bielefeld University, Germany,
goetze@math-uni.bielefeld.de.
}

\thanks{3) Research supported by the NSF grant DMS-2154001,
the GRF – SFB 1283/2 2021 – 317210226, and 
the Hausdorff Research Institute for Mathematics}


\begin{abstract}
An upper bound for the Kantorovich transport distance between probability
measures on multidimensional Euclidean spaces is given in terms of transport 
distances between one dimensional projections. This quantifies the
Cram\'er-Wold continuity theorem for the weak convergence of probability measures.
\end{abstract}

\maketitle
\markboth{Sergey G. Bobkov and Friedrich G\"otze}{Cram\'er-Wold theorem}

\def\theequation{\thesection.\arabic{equation}}
\def\E{{\mathbb E}}
\def\R{{\mathbb R}}
\def\C{{\mathbb C}}
\def\P{{\mathbb P}}
\def\Z{{\mathbb Z}}
\def\L{{\mathbb L}}
\def\T{{\mathbb T}}

\def\G{\Gamma}

\def\Ent{{\rm Ent}}
\def\var{{\rm Var}}
\def\Var{{\rm Var}}

\def\H{{\rm H}}
\def\Im{{\rm Im}}
\def\Tr{{\rm Tr}}
\def\s{{\mathfrak s}}

\def\k{{\kappa}}
\def\M{{\cal M}}
\def\Var{{\rm Var}}
\def\Ent{{\rm Ent}}
\def\O{{\rm Osc}_\mu}

\def\ep{\varepsilon}
\def\phi{\varphi}
\def\vp{\varphi}
\def\F{{\cal F}}

\def\be{\begin{equation}}
\def\en{\end{equation}}
\def\bee{\begin{eqnarray*}}
\def\ene{\end{eqnarray*}}

\thispagestyle{empty}

\vskip5mm
\section{{\bf Introduction}}
\setcounter{equation}{0}

\vskip2mm
\noindent
Given a sequence of random vectors $(X_n)_{n \geq 1}$ and a random vector $X$
with values in $\R^d$, the Cram\'er-Wold continuity theorem indicates that 
$X_n \Rightarrow X$ weakly in distribution, if and only if this convergences 
holds true for all one dimensional projections, i.e. if and only if
$$
\left<X_n,\theta\right> \Rightarrow \left<X,\theta\right> \quad {\rm as} \ n \rightarrow \infty
$$ 
on the real line for any $\theta \in \R^d$ (cf. \cite{C-W}, \cite{Bil}). 
Of  large interest is the problem of how one can quantify
this characterization by means of various distances responsible for the weak convergence.
Indeed, this could potentially reduce a number of high dimensional questions
to dimension one, perhaps under proper moment assumptions. Here we consider 
the problem with respect to the Kantorovich transport distance.

Let $X$ and $Y$ be random vectors in $\R^d$ with distributions
$\mu$ and $\nu$ having finite first absolute moments. 
The Kantorovich transport distance, also called the minimal distance between $\mu$
and $\nu$, is defined with respect to the Euclidean metric on $\R^d$ by
\be
W(X,Y) = W(\mu,\nu) = \inf\, \E\,|X' - Y'|= \inf \int_{\R^d}\int_{\R^d} |x-y|\,d\pi(x,y).
\en
Here the first infimum is taken over all pairs of random vectors $X',Y'$ with
distributions $\mu,\nu$, and the second one is running over all Borel probability
measures $\pi$ on $\R^d \times \R^d$ with marginals $\mu$ and $\nu$. 
By a well-known characterization of the convergence in $W$,
this metric metrizes the topology of weak convergence in the space of all Borel probability
measures $\mu$ on $\R^d$ with bounded $p$-th absolute moments for any fixed $p>1$ 
(cf. \cite{V}, Theorem 7.12).

According to the Kantorovich duality theorem (cf. \cite{D2}),
\begin{eqnarray}
W(X,Y) 
 & = &
\sup_{\|u\|_{\rm Lip} \leq 1}\, |\E u(X) - \E u(Y)| \nonumber \\
 & = &
\sup_{\|u\|_{\rm Lip} \leq 1}\Big|\int_{\R^d} u\,d\mu - \int_{\R^d} u\,d\nu\Big|,
\end{eqnarray}
where the supremum runs over all functions $u:\R^d \rightarrow \R$
with Lipschitz semi-norm $\|u\|_{\rm Lip} \leq 1$.
Since the functions of the form $u(x) = v(\left<x,\theta\right>)$ with $|\theta|=1$ and
$v:\R \rightarrow \R$ such that $\|v\|_{\rm Lip} \leq 1$
participate in the supremum (1.2), we have
\be
W(X,Y) \geq \sup_{|\theta| = 1}  W(X_\theta,Y_\theta)
\en
for the linear functionals
$$
X_\theta = \left<X,\theta\right>, \ Y_\theta = \left<Y,\theta\right>.
$$
In the modern literature, the supremum in (1.3) is often called the
max-scliced Wasserstein distance.
Recall that in dimension one, in big contrast to the multidimensional
situation, the Kantorovich distance has a simple description
$$
W(X_\theta,Y_\theta) = \int_{-\infty}^\infty |F_\theta(x) - G_\theta(x)|\,dx
$$
in terms of the distribution functions
$F_\theta(x) = \P\{X_\theta \leq x\}$, $G_\theta(x) = \P\{Y_\theta \leq x\}$.

Here we reverse the inequality (1.3) in a somewhat 
similar form under a $p$-th moment assumption.

\vskip2mm
{\bf Theorem 1.1.} {\sl Suppose that $(\E\,|X|^p)^{1/p} \leq b$ and 
$(\E\,|Y|^p)^{1/p} \leq b$ for some $p>1$ and $b \geq 0$. Then
\be
W(X,Y) \leq 18\,b^{1-\alpha} \sup_{|\theta| = 1}  W(X_\theta,Y_\theta)^\alpha,
\en
where
\be
\alpha = \frac{2}{d p^* + 2}, \quad p^* = \frac{p}{p-1}.
\en
}
\vskip3mm
Note that the distance $W$ is homogeneous with respect to $(X,Y)$, and so is
(1.4) when this inequality is written with an optimal value of $b$.

Letting $p \rightarrow \infty$ and assuming that $|X| \leq 1$ and $|Y| \leq 1$ a.s.,
we obtain a simpler relation
\be
W(X,Y) \leq 18\, \sup_{|\theta| = 1}  W(X_\theta,Y_\theta)^{\frac{2}{d+2}}.
\en
A similar relation with averaging over $\theta$ in place of the supremum in (1.3),
$$
W(X,Y) \, \leq \, c_d \, \Big(\int_{S^{d-1}} 
W(X_\theta,Y_\theta) \,d\sigma_{d-1}(\theta)\Big)^{\frac{1}{d+1}}
$$
was earlier obtained by Bonnotte  \cite{Bon}. Here $\sigma_{d-1}$ denotes the uniform
distribution on the unit sphere $S^{d-1}$ in $\R^d$, and $c = c_d$ is some
constant depending on the dimension $d$.

As another interesting case, suppose that $\E\,|X|^2 \leq d$ and $\E\,|Y|^2 \leq d$
(which holds for isotropic distributions). Then, using $d^{\frac{d}{2d+2}} \leq \sqrt{d}$
so that to simplify the factor $b^{1-\alpha}$, (1.4) yields
$$
W(X,Y) \leq 18\sqrt{d}\, \sup_{|\theta| = 1}  W(X_\theta,Y_\theta)^{\frac{1}{d+1}}.
$$

{\bf Remark 1.2.}
One may wonder whether or not it is possible to replace the Kantorovich distance
in Theorem 1.1 with other classical distances such as, for example, the multivariate
Kolmogorov distance
$$
\rho_d(X,Y) = \sup |\P\{X \in H\} - \P\{Y \in H\}|,
$$
where the supremum is running over all sets of the form
$H = (-\infty,x_1] \times \dots \times (-\infty,x_d]$ in $\R^d$. In connection with the
Cram\'er-Wold theorem, the computer tomography problems, and continuity properties 
of the Radon transform (which assigns to every probability distribution on $\R^d$ 
the collection of all its one-dimensional projections), this question and related inversion issues 
were discussed in 1990's in a series of works by Zinger, Klebanov and Khalfin,
cf. e.g. \cite{K-K}. However, as was demonstrated by Zaitsev \cite{Z}, the smallness of 
$\sup_{|\theta| = 1} \rho_1(X_\theta,Y_\theta)$ does not guarantee that 
$\rho_d(X,Y)$ will be small as well, even if the distributions of $X$ and $Y$ are compactly supported.
Moreover, here $\rho_1$ may be strengthened to the total variation distance.
The corresponding counter-example shows that $\rho_d$ is essentially stronger
than the Kantorovich metric $W$.

\vskip7mm
\section{{\bf Transport Distances and Convergence of Empirical Measures}}
\setcounter{equation}{0}

\vskip2mm
\noindent
It is not clear whether or not the exponent $\alpha = \alpha(p,d)$ defined in (1.5) is optimal in
the inequality (1.4), even if one can add an additional $(p,d)$-dependent multiplicative factor. 
In order to illustrate the strength of the inequality (1.6), we consider the following example involving 
empirical measures.

Let $X_1,\dots,X_n$ be a sample of size $n$ drawn from $\mu$, that is, independent random 
vectors in $\R^d$ with distribution $\mu$. They may be treated as independent copies of 
a random vector $X$. The associated empirical measures are defined by
\be
\mu_n = \frac{1}{n} \sum_{k=1}^n \delta_{X_k},
\en
where $\delta_x$ denotes a delta-measure at the point $x$. Correspondingly,
their linear projections represent one dimensional empirical measures
\be
\mu_{n,\theta} = \frac{1}{n} \sum_{k=1}^n \delta_{\left<X_k,\theta\right>}.
\en

For simplicity, let us restrict ourselves to the distributions 
$\mu$ supported on the unit ball $B_1$ in $\R^d$. 
By concavity of
$x\to x^{\alpha}$ in $x \geq 0$ for $\alpha \in (0,1]$, it follows from (1.6) that
\be
\E\,W(\mu_n,\mu) \leq 18\, \Big[ \E \sup_{|\theta| = 1} 
W(\mu_{n,\theta},\mu_\theta)\Big]^\alpha
\en
with $\alpha = \frac{2}{d+2}$, 
In the case $d \geq 3$, it is known (cf. e.g. \cite{D1}, \cite{B-L2}) that $\E\,W(\mu_n,\mu)$ 
is of order at most $c_d n^{-1/d}$, and this rate cannot be improved as $n \rightarrow \infty$ 
for the uniform distribution. However, if $d=1$ and $\mu$ is compactly supported, 
the rate is of the standard order
$\E\,W(\mu_n,\mu) \sim \frac{1}{\sqrt{n}}$ up to $\mu$-dependent factors 
(for a general two-sided bound we refer to \cite{B-L1}, Theorem 3.5).
Hence, we also have $\E\,W(\mu_{n,\theta},\mu_\theta) \sim \frac{1}{\sqrt{n}}$
for every fixed $\theta$. In order to determine the upper bound  for the transport distance 
via the bound (2.3), the following result can be used.

\vskip5mm
{\bf Theorem 2.1.} {\sl If $\mu$ is supported on the unit ball $B_1$
in $\R^d$, then
\be
\E \sup_{|\theta| = 1} W(\mu_{n,\theta},\mu_\theta) \leq \frac{c_d}{\sqrt{n}}
\en
with some constants $c_d>0$ depending on $d$ only.
}

\vskip5mm
Applying this bound in (2.3) and using a lower bound with rate $c_d n^{-1/d}$ 
for the left-hand side with large $n$, we may conclude that necessarily $\alpha \leq 2/d$.
Thus, the exponent $\alpha = \frac{2}{d+2}$ is asymptotically optimal
for the growing dimension $d$.

For reader's convenience, at the end of the note we include a simple chaining argument leading to (2.4) with
$c_d = c\sqrt{d}$, where $c>0$ is an absolute constant. However, this statement is not new --
inequalities for one-dimensional projections of empirical measures such as (2.4) have been 
the subject of many recent investigations, cf. e.g. \cite{NW-R}, \cite{N-G-S-K}, \cite{M-B-W}.
After this paper was submitted, we also learned about the preprint by Boedihardjo \cite{Boe},
where the upper bound (2.4) was derived with a constant independent of $d$.

\vskip7mm
\section{{\bf Reduction to Compactly Supported Lipschitz Functions}}
\setcounter{equation}{0}

\vskip2mm
\noindent
We need some preparation for the proof of Theorem 1.1. The argument is based on truncation, 
smoothing, the Plancherel theorem, together with the Kantorovich duality theorem (1.2).

Let $U_r$ ($r>0$) denote the collection of all functions $u:\R^d \rightarrow \R$ with 
$\|u\|_{\rm Lip} \leq 1$, $u(0) = 0$, which are supported on the
Euclidean ball $B_r = \{x \in \R^d: |x| \leq r\}$. As a first step, we show that
the supremum in (1.2) may be restricted to the set $U_r$ at the expense of a small error
for large values of the parameter $r$ under a $p$-th moment assumption. Define
\bee
W^{(r)}(X,Y) = \sup_{u \in U_r}\, |\E u(X) - \E u(Y)| = \sup_{u \in U_r}
\Big|\int_{\R^d} u\,d\mu - \int_{\R^d} u\,d\nu\Big|,
\ene
assuming that the random vectors $X$ and $Y$ have distributions $\mu$ and $\nu$.

\vskip5mm
{\bf Lemma 3.1.} {\sl Let $(\E\,|X|^p)^{1/p} \leq b$ and $(\E\,|Y|^p)^{1/p} \leq b$ for 
some $p>1$. For any $r>0$,
$$
W(X,Y) \leq 3\, W^{(r)}(X,Y) + 4b\, \Big(\frac{2b}{r}\Big)^{p-1}.
$$
}

{\bf Proof.} Let $u$ be a Lipschitz function participating in the supremum (1.2) with 
$u(0) = 0$ (without loss of generality). The latter ensures that $|u(x)| \leq |x|$ for all $x \in \R^d$.

Define the function $u_r = u \psi_r$, where using the notation $a^+ = \max(a,0)$,
$$
\psi_r(x) = \Big(1 - \frac{2}{r}\, {\rm dist}(B_{r/2},x)\Big)^+, \quad x \in \R^d.
$$
Clearly, $0 \leq \psi_r \leq 1$ and $\|\psi_r\|_{\rm Lip} \leq \frac{2}{r}$, as the distance
function 
$$
x \rightarrow {\rm dist}(B,x) = \inf\{|x-y|: y \in B\}
$$
has a Lipschitz semi-norm at most one for any non-empty set $B$ in $\R^d$.

By the definition, $|u_r(x)| \leq |u(x)|$ for all $x \in \R^d$, and
\be
u_r(x) = u(x) \ {\rm \ for} \ |x| \leq r/2, \qquad u_r(x) = 0 \ {\rm \ for} \ |x| \geq r.
\en
Writing
$$
u_r(x) - u_r(y) = (u(x) - u(y)) \,\psi_r(x) + u(y)\,(\psi_r(x) - \psi_r(y)),
$$
it follows that, for all $x \in \R^d$ and $y \in B_r$,
$$
|u_r(x) - u_r(y)| \leq |u(x) - u(y)| + |y|\, \frac{2}{r}\, |x-y| \leq 3\,|x-y|.
$$
A similar final inequality holds true for $x \in B_r$ and $y \in \R^d$.
In addition, $u_r(x) - u_r(y) = 0$ in the case $x,y \notin B_r$.
Therefore, $\|u_r\|_{\rm Lip} \leq 3$.

Next, by (3.1), 
\bee
\Big|\int_{\R^d} u\,d\mu - \int_{\R^d} u_r\,d\mu\Big| 
 & \leq &
\int_{|x| > r/2} |u(x)|\,d\mu(x) + \int_{|x| > r/2} |u_r(x)|\,d\mu(x) \\
 & \leq &
2\int_{|x| > r/2} |x|\,d\mu(x) \\
 & = & 
2\,\E\,|X|\,1_{\{|X| > r/2\}} \, \leq \, 2\,\frac{\E\,|X|^p}{(r/2)^{p-1}}\, \leq \,
2b^p\, \Big(\frac{2}{r}\Big)^{p-1}.
\ene
With a similar inequality for the measure $\nu$, we obtain that
$$
\Big|\int_{\R^d} u\,d(\mu - \nu) - \int_{\R^d} u_r\,d(\mu - \nu)\Big|\, \leq \,
4b^p\, \Big(\frac{2}{r}\Big)^{p-1},
$$
implying
$$
\Big|\int_{\R^d} u\,d(\mu - \nu)\Big| \leq \Big|\int_{\R^d} u_r\,d(\mu - \nu)\Big| +
4b^p\, \Big(\frac{2}{r}\Big)^{p-1}.
$$
Since $\|u_r\|_{\rm Lip} \leq 3$, we have $\frac{1}{3}\,u_r \in U_r$, so that
the last integral does not exceed $3\, W^{(r)}(X,Y)$ in absolute value. Thus,
$$
\Big|\int_{\R^d} u\,d(\mu - \nu)\Big|\leq 3\, W^{(r)}(X,Y) + 4b^p\, \Big(\frac{2}{r}\Big)^{p-1}.
$$
It remains to take the supremum on the left-hand side over all functions 
$u:\R^d \rightarrow \R$ such that $\|u\|_{\rm Lip} \leq 1$ and $u(0) = 0$.
\qed

\vskip7mm
\section{{\bf Fourier Transforms}}
\setcounter{equation}{0}

\vskip2mm
\noindent
Any integrable compactly supported function $u$ on $\R^d$ has 
a well-defined Fourier transform
\be
\widehat u(t) = \int_{\R^d} e^{i \left<t,x\right>} u(x)\,dx, \quad t \in \R^d,
\en
which represents a $C^\infty$-smooth function. 
Towards the proof of Theorem 1.1 let us state now the following integrability property.

\vskip5mm
{\bf Lemma 4.1.} {\sl For any function $u:\R^d \rightarrow \R$ which is supported 
on the ball $B_r$ and has a Lipschitz semi-norm $\|u\|_{\rm Lip} \leq 1$,
\be
\int_{\R^d} |\widehat u(t)|^2\, |t|^2\,dt \, \leq \, \omega_d\, (2\pi r)^d.
\en
}
\vskip2mm
In particular, this inequality holds true for any function $u$ in $U_r$.
Here and elsewhere $\omega_d$
stands for the $d$-dimensional volume of the unit ball $B_1$. 

\vskip5mm
{\bf Proof.} First assume that $\widehat u(t)$ decays sufficiently fast at infinity,
namely $\widehat u(t) = O(1/|t|^p)$ as $|t| \rightarrow \infty$ for any $p>0$.
Then (4.1) may be inverted in the form of the Fourier transform
$$
u(x) = \frac{1}{(2\pi)^d} \int_{\R^d} e^{-i \left<t,x\right>}\, \widehat u(t)\,dt.
$$
In particular, $u$ is $C^\infty$-smooth on $\R^d$. Moreover, this equality may be
differentiated along every coordinate $x_k$ to represent the corresponding partial
derivatives as
$$
\partial_{x_k} u(x) = -\frac{i}{(2\pi)^d} \int_{\R^d} e^{-i \left<t,x\right>}\,t_k\, \widehat u(t)\,dt,
\quad k = 1,\dots,d.
$$
Hence, by the Plancherel theorem,
$$
\int_{\R^d} t_k^2\, |\widehat u(t)|^2\,dt = (2\pi)^d \int_{\R^d} (\partial_{x_k} u(x))^2\,dx.
$$
Summing over all $k \leq d$ and using $|\nabla u(x)| \leq 1$ for $x \in B_r$ and
$\nabla u(x) = 0$ for $|x| > r$, we get
$$
\int_{\R^d} |t|^2\,|\widehat u(t)|^2\, dt = (2\pi)^d
\int_{\R^d} |\nabla u(x)|^2\,dx \leq (2\pi)^d \cdot \omega_d r^d,
$$
which is the desired inequality (4.2).

In the general case, a smoothing argument can be used. By the well-known theorem
of Ingham \cite{I}, there exists a probability density $w$ 
on $\R^d$ which is supported on the unit ball $B_1$ and has characteristic function 
$\widehat w(t)$ satisfying $\widehat w(t) = O(1/|t|^p)$ as $|t| \rightarrow \infty$,
for any fixed $p>0$.
Given $\ep>0$, the probability density $w_\ep(x) = \ep^{-d} w(x/\ep)$ is supported 
on the ball $B_\ep$ and has characteristic function $\widehat w_\ep(t) = \widehat w(\ep t)$. 
Consider the convolution
$$
u_\ep(x) = (u * w_\ep)(x) = \int_{\R^d} u(x-y)\,w_\ep(y)\,dy, \quad x \in \R^d.
$$
This function is supported on $B_{r+\ep}$ and has a Lipschitz semi-norm
$\|u_\ep\|_{\rm Lip} \leq \|u\|_{\rm Lip} \leq 1$. By the Lipschitz property,
$|u(x)| \leq |u(0)| + r$ for any $x \in B_r$, implying that
$$
\sup_{t \in \R^d}
|\widehat u(t)| \leq \int_{B_r} |u(x)|\,dx \leq (|u(0)| + r)\,\omega_d r^d < \infty.
$$
Hence, the Fourier transform of $u_\ep$ satisfies
$\widehat u_\ep(t) = \widehat u(t) \widehat w(\ep t) = O(1/|t|^p)$ as $|t| \rightarrow \infty$.
Thus, one may apply the previous step to the function $u_\ep$ which gives
$$
\int_{\R^d} |\widehat u(t)|^2\, |\widehat w(\ep t)|^2\,|t|^2\,dt \, \leq \, 
\omega_d\, (2\pi\, (r+\ep))^d.
$$
It remains to send $\ep \rightarrow 0$ in this inequality and apply Fatou's lemma together with
$\widehat w(\ep t) \rightarrow 1$ as $\ep \rightarrow 0$.
\qed

\vskip4mm
Next, let us connect the Kantorovich distance with the multivariate characteristic functions
$$
f(t) = \E\,e^{i\left<t,X\right>}, \quad g(t) = \E\,e^{i\left<t,Y\right>} \qquad (t \in \R^d).
$$

\vskip5mm
{\bf Lemma 4.2.} {\sl Given random vectors $X,Y$ in $\R^d$ with characteristic
functions $f,g$ and finite first absolute moments, we have, for any $t \in \R^d$,
$$
|f(t) - g(t)| \leq 2|t|\,\sup_{|\theta| = 1}  W(X_\theta,Y_\theta).
$$
}

{\bf Proof.} In dimension one, using the property that the function $u_t(x) = \frac{1}{t}\,\cos^{itx}$
with parameter $t \neq 0$, has a Lipschitz semi-norm at most 1, it follows from (1.2) that
$$
|{\rm Re}(f(t)) - {\rm Re}(g(t))| \leq |t|\,W(X,Y).
$$
By a similar argument,
$$
|{\rm Im}(f(t)) - {\rm Im}(g(t))| \leq |t|\,W(X,Y),
$$
so that
$$
|f(t) - g(t)| \leq 2|t|\,W(X,Y).
$$

In dimension $d$, just note that, for any $\theta \in \R^d$, the functions
$r \rightarrow f(r\theta)$ and $r \rightarrow g(r\theta)$
represent the characteristic functions of $X_\theta$ of $Y_\theta$.
\qed

\vskip7mm
\section{{\bf Proof of Theorem 1.1}}
\setcounter{equation}{0}

\noindent
With $b = \max(\|X\|_p,\|Y\|_p)$, where
$$
\|X\|_p = (\E\,|X|^p)^{1/p}, \quad \|Y\|_p = (\E\,|Y|^p)^{1/p},
$$
the inequality (1.4) is homogeneous with respect to $(X,Y)$, so one may assume that $b=1$.
As a consequence, 
$$
W(X,Y) \leq \E\,|X| + \E\,|Y| \leq 2.
$$

Let $\eta$ be a random vector with uniform distribution in the
ball $B_1$, that is, with density $w(x) = \frac{1}{\omega_d}\,1_{B_1}(x)$,
and let $h(t)$ denote its characteristic function. Consider the random vectors
$$
X(\ep) = X + \ep \eta, \quad Y(\ep) = Y + \ep \eta \quad (\ep > 0),
$$
assuming that $\eta$ is independent of $X$ and $Y$. Then, by the definition (1.1),
or by (1.2),
\be
W(X,Y) \leq W(X(\ep),Y(\ep)) + 2\ep.
\en
Indeed, given a function $u$ on $\R^d$ with $\|u\|_{\rm Lip} \leq 1$, we have
\bee
\E\, u(X(\ep)) - \E\, u(Y(\ep))
 & = &
\E\, u(X + \ep \eta) - \E\, u(Y + \ep \eta) \\
 & \geq &
\E\, \big(u(X) - \ep |\eta|\big) - \E\, \big(u(Y) + \ep |\eta|\big) \\
 & \geq &
\E\,u(X) - \E\,u(Y) - 2\ep.
\ene
Taking the supremum of both sides over all Lipschitz $u$ and applying (1.2), we arrive at (5.1).

On the other hand, using
$$
\|X(\ep)\|_p \leq 1+\ep, \quad \|Y(\ep)\|_p \leq 1+\ep,
$$
one may apply Lemma 3.1, which gives that, for any $r>0$,
$$
W(X(\ep),Y(\ep)) \leq 3\, W^{(r)}(X(\ep),Y(\ep)) + 
4(1+\ep)\, \Big(\frac{2(1+\ep)}{r}\Big)^{p-1}.
$$
Therefore, by (5.1),
\be
W(X,Y) \leq 3\, W^{(r)}(X(\ep),Y(\ep)) + 
4(1+\ep)\, \Big(\frac{2(1+\ep)}{r}\Big)^{p-1} + 2\ep.
\en

In order to estimate the first term on the right-hand side, first note that 
the distributions of $X(\ep)$ and $Y(\ep)$ represent convolutions of the 
distributions of $X$ and $Y$ with a uniform distribution on the ball $B_\ep$.
Hence, these random vectors have densities which we denote by $p_\ep$ and $q_\ep$ 
respectively. They have respective characteristic functions
$$
f_\ep(t) = f(t) h(\ep t), \quad g_\ep(t) = g(t) h(\ep t) \qquad (t \in \R^d),
$$
where $f$ and $g$ denote the characteristic functions of $X$ and $Y$.
As the function $h(t)$ is square integrable, while $|f(t)| \leq 1$ and
$|g(t)| \leq 1$, the functions $f_\ep$ and $g_\ep$ are square integrable, so that
the densities $p_\ep$ and $q_\ep$ are square integrable as well, by the Plancherel 
theorem. Thus, given a function $u$ in $U_r$, one may write
$$
\E\, u(X(\ep)) - \E\, u(Y(\ep)) = \int_{\R^d} u(x)\,(p_\ep(x) - q_\ep(x))\,dx.
$$

We are in position to apply the Plancherel theorem once more and 
rewrite the last integral as
$$
\frac{1}{(2\pi)^d} \int_{\R^d} \widehat u(t)\,(\bar f(t) - \bar g(t))\,h(\ep t)\,dt.
$$
Thanks to Lemma 4.2 we then have
\be
|\E\, u(X(\ep)) - \E\, u(Y(\ep))| \leq 
\frac{2M}{(2\pi)^d} \int_{\R^d} |\widehat u(t)|\,|t|\,|h(\ep t)|\,dt,
\en
where 
$$
M = \sup_{|\theta| = 1}  W(X_\theta,Y_\theta).
$$ 
Moreover, using
\bee
\int_{\R^d} |h(\ep t)|^2\,dt 
 & = &
\ep^{-d} \int_{\R^d} |h(t)|^2\,dt \\
 & = &
(2\pi)^d\, \ep^{-d} \int_{\R^d} w(x)^2\,dx \, = \, (2\pi)^d\, \ep^{-d} \omega_d^{-1}
\ene
and applying Cauchy's inequality in (5.3) together with Lemma 4.1, we obtain that
$$
|\E\, u(X(\ep)) - \E\, u(Y(\ep))| \leq 2M \Big(\frac{r}{\ep}\Big)^{\frac{d}{2}}.
$$
Taking the supremum over all $u \in U(r)$ on the left-hand side leads
to the similar bound for $W^{(r)}(X(\ep),Y(\ep))$, and using this in (5.2) we are led to
$$
W(X,Y) \, \leq \, 6M\, \Big(\frac{r}{\ep}\Big)^{\frac{d}{2}} + 
4(1+\ep)\, \Big(\frac{2(1+\ep)}{r}\Big)^{p-1} + 2\ep.
$$

To simplify optimization over free variables $r>0$ and $\ep>0$, let us assume that
$2(1+\ep) \leq c$ for a constant $c>2$ (to be chosen later on). Then we have
$$
W(X,Y) \, \leq \, 6M\, \Big(\frac{r}{\ep}\Big)^{\frac{d}{2}} + 
2c\, \Big(\frac{c}{r}\Big)^{p-1} + 2\ep.
$$
Let us then replace $r$ with $cs$ and $\ep$ with $c\delta$ in the above inequality to get 
$$
W(X,Y) \, \leq \, 6M\, \Big(\frac{s}{\delta}\Big)^{\frac{d}{2}} + 
2c\, \Big(\frac{1}{s}\Big)^{p-1} + 2c\delta.
$$
Here, equalizing the terms $M\,(\frac{s}{\delta})^{\frac{d}{2}}$ and $s^{-(p-1)}$,
we find the unique value of $s$ for which the above yields
\be
W(X,Y) \, \leq \, (6+2c)\, A \delta^{- \beta} + 2c\delta,
\en
where
$$
\quad \beta = \frac{(p-1)\,\frac{d}{2}}{p-1 + \frac{d}{2}}, \quad
A = M^{\frac{p-1}{p-1 + \frac{d}{2}}}.
$$
The choice $\delta = A^{\frac{1}{\beta + 1}}$ in (5.4) leads to
$$
W(X,Y) \, \leq \, (6+4c)\, A^{\frac{1}{\beta + 1}},
$$
provided that $2(1+c\delta) \leq c$. If we require that $\delta \leq \frac{1}{6}$,
the latter condition is satisfied for $c=3$, and we obtain that
\be
W(X,Y) \, \leq \, 18\, A^{\frac{1}{\beta + 1}},
\en
provided that $A^{\frac{1}{\beta + 1}} \leq \frac{1}{6}$.
In the other case, the right-hand side in (5.5) is greater than 2, so this 
inequality is fulfilled automatically due to the property $W(X,Y) \leq 2$.

It remains to note that
$$
A^{\frac{1}{\beta + 1}} = M^{\frac{2}{d p^* + 2}}, \quad
p^* = \frac{p}{p-1}.
$$
\qed

\vskip7mm
\section{{\bf Proof of Theorem 2.1}}
\setcounter{equation}{0}

\noindent
Let $U$ denote  for the space of all functions $u:[-1,1] \rightarrow \R$ with 
$\|u\|_{\rm Lip} \leq 1$, such that $u(0) = 0$. We equip $U$ with the uniform distance 
$$
\|u-v\|_\infty = \max\{|u(x) - v(x)|: |x|,|y|\leq 1\},
$$
which turns this set into a compact space, by the Arzel\'a-Ascoli theorem.

The empirical measures $\mu_n$ and  $\mu_{n,\theta}$ defined in (2.1)-(2.2)
are random with mean $\mu$ so that
$$
\int_{\R} u\,d\mu_{n,\theta} = \frac{1}{n} \sum_{k=1}^n u(\left<X_k,\theta\right>), \quad
\int_{\R} u\,d\mu_\theta = \E\, u(\left<X,\theta\right>).
$$
According to the Kantorovich duality theorem (1.2),
\be
\sup_{|\theta| = 1} W(\mu_{n,\theta},\mu_\theta) \, = \,
\frac{1}{n}\, \sup_{\theta \in S^{d-1}} \, \sup_{u \in U}\,\bigg|
\sum_{k=1}^n \big(u(\left<X_k,\theta\right>) - \E\, u(\left<X_k,\theta\right>)\big)\bigg|,
\en
where $S^{d-1} = \{\theta \in \R^d: |\theta|=1\}$ denotes the unit sphere in $\R^d$. 

In order to bound the expectation of the right-hand side in (6.1), one may use
chaining arguments, in particular, a well-known theorem by Dudley 
which says the following (for a proof, let us refer to \cite{T1}, \cite{T2}). 
Given a random variable $\xi$, define its Orlicz $\psi_2$-norm
$$
\|\xi\|_{\psi_2} = \inf\big\{\lambda>0: \E\, e^{\xi^2/\lambda^2} \leq 2\big\}.
$$
Suppose that $\xi(t)$ is a mean zero random process defined on some compact
metric space $(T,\rho)$, which satisfies the Lipschitz property
\be
\|\xi(t) - \xi(s)\|_{\psi_2} \leq \Lambda \rho(t,s), \quad t,s \in T,
\en
with some $\Lambda>0$. Then with some absolute constant $K$ we have
\be
\E\, \sup_{t \in T}\, \xi(t) \, \leq \, K\Lambda \int_0^D \sqrt{\log N(\ep)}\ d\ep,
\en
where $D = \max\{\rho(t,s): t,s \in T\}$ is the diameter and 
$N(\ep) = N(T,\rho,\ep)$ is the minimal number of closed balls 
in $T$ of radius $\ep$ needed to cover the space (recall that $\log N(\ep)$ is
called the $\ep$-entropy of $(T,\rho)$).

In view of (6.1), it is natural to consider the random process
$$
\xi(t) =  \xi(u,\theta) = \frac{1}{\sqrt{n}} \sum_{k=1}^n 
\big(u(\left<X_k,\theta\right>) - \E\, u(\left<X_k,\theta\right>)\big),
\quad t = (u,\theta) \in T = U \times S^{d-1}.
$$
We equip $T$ with the metric
$$
\rho(t,s) = \|u - v\|_\infty + |\theta - \theta'|, \quad
t = (u,\theta), s = (v,\theta') \in T.
$$
Endowed with this metric $T$ will be  a compact space of diameter $D=4$.

Now, given two points $t = (u,\theta)$, $s = (v,\theta')$ in $T$, one may write
$$
\xi(t) - \xi(s) = \frac{1}{\sqrt{n}}\, \sum_{k=1}^n (\eta_k - \E \eta_k),
$$
where
$$
\eta_k = u(\left<X_k,\theta\right>) - v(\left<X_k,\theta'\right>).
$$
These random variables are independent and bounded. Indeed, writing
$$
\eta_k = \big(u(\left<X_k,\theta\right>) - v(\left<X_k,\theta\right>)\big) +
\big(v(\left<X_k,\theta\right>) - v(\left<X_k,\theta'\right>)\big)
$$
and using the Lipschitz property of $v$ together with the assumption
$|X_k| \leq 1$ a.s., we have 
$$
|\eta_k| \leq \rho(t,s) \ {\rm a.s.} 
$$
Recall that, by the well-known Hoeffding’s lemma,
for any random variable $\eta$ such that $|\eta| \leq r$ a.s.,
$$
\E\, e^{z (\eta - \E \eta)} \leq e^{r^2 z^2/2} \quad {\rm for \ all} \ z \in \R.
$$ 
Hence, this holds for all $\eta = \eta_k$ with $r = \rho(t,s)$, and thus
$$
\E\, e^{z (\xi(t) - \xi(s))} \leq e^{\rho^2 z^2/2}, \quad \rho = \rho(t,s).
$$
Integrating this inequality with respect to $z$ over the Gaussian
measure with mean zero and variance $\sigma^2$ ($0 < \sigma < 1/\rho$),
we get
$$
\E\, e^{\sigma^2 (\xi(t) - \xi(s))^2/2} \leq \frac{1}{\sqrt{1 - (\sigma \rho)^2}}
$$ 
which implies
$$
\|\xi(t) - \xi(s)\|_{\psi_2} \leq \rho(t,s) \sqrt{8/3}.
$$

Thus, the Lipschitz condition (6.2) is fulfilled with an absolute constant $\Lambda$,
and we may apply the Dudley's bound (6.3). By the definition of the metric $\rho$
in $T$, any closed ball of radius $2\ep$ in this space contains the product
of a closed ball in $U$ and a closed ball in $S^{d-1}$, both of radius $\ep$.
Hence, the corresponding $\ep$-entropies are connected by the relation
\be
N(T,\rho,2\ep) \leq N(U,\ep)\, N(S^{d-1},\ep).
\en
It is well-known (cf. \cite{K-T}) that
$$
\log N(U,\ep) \leq \frac{c}{\ep}, \quad N(S^{d-1},\ep) \leq 
\Big(\frac{c}{\ep}\Big)^d, \quad 0 < \ep \leq 2,
$$
with some absolute constant $c>0$. Using these bounds in (6.4) and then
in (6.3), we conclude that the expectation of both sides in (6.1) does not
exceed a multiple of $\sqrt{d/n}$. This gives the desired relation (2.4).
\qed

\vskip5mm
{\bf Acknowledgement.} We would like to thank the referee for a careful reading
and useful comments. We also thank Andrey Zaitsev and March Boedihardjo for providing
us missing references in the original arxiv version of this paper.

\end{document}